\documentclass{commatDV}

\usepackage{DLde}
\usepackage{graphicx}

\newcommand\mkA{\ensuremath{\mathfrak a}}
\DeclareMathOperator\END{End}
\DeclareMathOperator\im{Im}

\title{Transitive irreducible Lie superalgebras of vector fields}

\author[Arkady Onishchik]{\fbox{Arkady Onishchik}}

\affiliation{N/A}

\abstract{Let $\mathfrak{d}$  be the Lie superalgebra of superderivations of the sheaf of sections of the exterior algebra of the homogeneous vector bundle $E$ over the flag variety $G/P$, where $G$ is a~simple finite-dimensional complex Lie group and $P$ its parabolic subgroup. Then, $\mathfrak{d}$ is transitive and irreducible whenever $E$  is defined by an irreducible $P$-module $V$ such that the highest weight of $V^*$ is dominant. Moreover,  $\mathfrak{d}$ is simple; it is isomorphic to the Lie superalgebra of vector fields on the superpoint, i.e., on a~$0|n$-dimensional supervariety.
}



\VOLUME{30}
\NUMBER{3}
\firstpage{25}
\DOI{https://doi.org/10.46298/cm.10456}

\begin{paper}

\section*{Preface of the editor}\label{S:0}
The manuscript of this paper was deposited in VINITI 12.06.86 No 4329-B which is inaccessible. There still appear papers with references  to some parts of this inaccessible deposition, see, e.g.,  \cite{IO}, \cite{OVs} which contains references to this preprint and its continuation \cite{On} preprinted in the proceedings of  the  ``Seminar on Supersymmetries'', see \url{http://staff.math.su.se/mleites/sos.html}. So I decided to make it available, together with \cite{On}, by translating these texts and updating the references; the ones I added are endowed with an asterisk. The abstract, footnotes, and comments are due to me. 

For a~comprehensive description of simple Lie superalgebras of vector fields over algebraically closed fields of any characteristic, see \cite{BGLLS}. \textit{D. Leites}

\section{Basics. Introduction}\label{S:1}

Recall that a~ filtered Lie superalgebra 
\[
\mkA =\mkA _{(-d)}\supset \dots \supset \mkA _{(-1)}\supset \mkA _{(0)}\supset \mkA _{(1)}\supset \ldots , \text{~~where $d\in\Nee:=\{1,2,\dots\}$}
\]
is said to be \textit{transitive}, see \cite{Dr}, if  
\begin{equation}
\label{11}
\mkA _{(p+1)}=\{x\in \mkA _{(p)}\mid [x,\mkA ]\subset \mkA _{(p)}\} \text{~~for all $p\geq 0$}.
\end{equation}

A $\Zee$-graded Lie superalgebra
$\fb=\oplus_{i \geq -d}\fb_i$
is said to be \textit{transitive} if for all $p\geq 0$ we have
\begin{equation}
\label{12}
\{x\in \fb_p \mid [x,\fb_{-}]=0\}=0, \text{~~where $\fb_{-}:=\oplus _{p <0}\fb_p$}.
\end{equation}
As is well-known, see \cite{K}, \cite{Sch}, 
condition ~\eqref{11} is satisfied if and only if the graded Lie algebra $\fb:=\gr \mkA $ satisfies conditions ~\eqref{12}.

We say that a graded Lie superalgebra  $\fb$ is \textit{irreducible} if the \textbf{adjoint} representation of $\fb _0$ in $\fb_{-1}$ is irreducible. (The term \textit{irreducible} is usually used to denote the image of the (Lie) algebra or superalgebra $\fb$ in any irreducible representation. \textit{D.L.}) 

Simple, and even primitive\footnote{Let a~Lie algebra $\cL$ contain a~subalgebra $\cL_0$ which does not contain any nonzero ideal of $\cL$. Let
\[
\cL_1:=\{t\in \cL_0\mid [t,\cL] \subset \cL_0\}.
\]
Then, $\cL$ is called \textit{primitive} (see \cite{O}) if 
\begin{equation}\label{prim}
\text{
$\cL_0$ is a~maximal subalgebra of $\cL$ and $\cL_1\neq \{0\}$.} 
\end{equation}

Unlike primitive Lie algebras (both finite-dimensional and $\Zee$-graded of polynomial growth) which do not differ much from the simple Lie algebras, the classification of primitive Lie \textbf{super}algebras over $\Cee$ announced by Kac, see \cite{K}, is a~wild problem, as shown in \cite{ALSh}. Compare with \cite{CaKa}.},  Lie (super)algebras have a~filtration in which the associated graded Lie (super)algebra is both transitive and irreducible.  For the classification of \textbf{simple} finite-dimensional Lie superalgebras over $\Cee$, see  \cite{K}, \cite{Sch}. (Actually, no complete proof of the classification of simple finite-dimensional  Lie superalgebras over $\Cee$ is known to this day; e.g., the proof of completeness of the list of examples of known deformations with odd parameters is not published. Such deformations were not even mentioned in  \cite{K}, \cite{Sch}. \textit{D.L.})

Hereafter I assume that the $\Zee$-grading of every Lie superalgebra considered is compatible with parity.

In this paper, I show that such Lie superalgebras naturally appear in the study of homogeneous vector bundles over complex homogeneous flag varieties $G/P$, where $P$ is a~parabolic subgroup of a~complex semisimple Lie group $G$. Namely, 
to each homogeneous vector bundle $E \tto M$ there is associated a~\textit{split} complex supermanifold $(M,\Lambda (\mathcal E))$, where $\mathcal E $ is the sheaf of holomorphic sections of $E$, and $\Lambda(\mathcal E)$ is the Grassmann algebra of $\mathcal E$. 

Here, I prove that the graded Lie superalgebra $\fd =\Gamma (M,\Der \Lambda(\mathcal E))$ of vector fields on this supermanifold is transitive and irreducible if $E$ is defined by an irreducible representation $\varphi $ of $P$ satisfying certain natural requirements. The main point here is that the highest weight of the representation $\varphi ^*$ should be \textit{dominant}. This guarantees the existence of sufficiently many holomorphic sections of ``the odd tangent bundle'' $E^*$ over the supermanifold $(M,\Lambda(\mathcal E))$.

The Bott--Borel--Weil theorem, along with other standard methods, enables one  to explicitly compute the Lie superalgebra $\fd$ in many cases.

In this paper, this computation is carried out in the cases where $E$ is the cotangent bundle over $M$, and $G$ is simple. 

\section{Vector fields on supermanifolds}\label{S:2}

The term ``supermanifold'' is meant in the sense of Berezin-Leites,  in the complex--analytic situation, see \cite{L}, \cite{MaG}, \cite{Va}. Let us recall the definition.

Let $(n, m)$ be a~pair of non-negative integers. Let a~\textit{model superspace} of dimension $n|m$ be the ringed space $(\Cee^n, \widetilde\cO)$, where $\widetilde\cO:=\Lambda_{\cO}[\xi]$ for $\xi=(\xi_1,\dots, \xi_m)$, and $\cO$ is the sheaf of germs of holomorphic functions on $\Cee^n$. A topological space $M$ equipped with a~sheaf of supercommutative superalgebras $\widetilde\cO$,  which as a~ ringed space is locally isomorphic to  a~model superspace of dimension $n|m$ is said to be a(n almost) \textit{complex}  \textit{supermanifold} and $\widetilde\cO$ is its \textit{structure sheaf}. (For details, e.g., which of almost complex supermanifolds are complex, and for phenomena indigenous \textbf{not only} to the ``super'' world --- \textit{real--complex (super)manifolds}, see \cite{BGLS}. \textit{D.L.})

Let us introduce \textit{parity} in $\widetilde\cO$ by declaring all of the $\xi_i$ to be odd; so $\widetilde\cO=\widetilde\cO_\ev\oplus \widetilde\cO_\od$ is the sum of its homogeneous components even and odd.

Let $\cI\subset \widetilde\cO$ be the subsheaf of ideals generated by subsheaf $\widetilde\cO_\od$ and let  $\widetilde\cO_{\rd}:= \widetilde\cO/\cI$. Clearly, $(M, \widetilde\cO_{\rd})$ is a(n almost) complex-analytic manifold.

The simplest supermanifolds are the \textit{split} ones which are described as follows. Let $(M, \cO)$ be a(n almost) complex-analytic manifold of dimension $n$ and $\cF$ a~locally free analytic sheaf of rank $m$ on $M$. Set $\widetilde\cO:=\Lambda_{\cO}(\cF)$. Clearly, $(M, \widetilde\cO)$ is a~supermanifold, and $\widetilde\cO_{\rd}$ is naturally isomorphic to the subsheaf $\cO\subset\widetilde\cO$. In what follows I will identify  $\widetilde\cO_{\rd}$ with $\cO$ in this situation and will briefly write  $\widetilde\cO=\Lambda_{\cO}(\cF)$. The structure sheaf of the split supermanifold is endowed with a~$\Zee$-grading
\[
\widetilde\cO=\bigoplus\limits _{0\leq p\leq m}\widetilde\cO_p, \text{~~where $\widetilde\cO_p=\Lambda ^p_{\cO}(\cF)$}.
\]
Observe that \textbf{every $C^\infty$ of analytic supermanifold is locally split}.\footnote{This is not true for algebraic supervarieties and supeschemes, as follows from \cite{AD}: over a~contractible paracompact real set all vector bundles are trivial, but this is not necessarily true over other fields and non-affine schemes. The first example of non-split supervariety is due to P.~Green (\cite{Gr}), see also \cite{Pa} submitted only 2 month later than  \cite{Gr} and expounded in \cite[Ch.4, \S4, Sections 6--9]{Ber1} reproduced in \cite[Ch.3, Th.2, p.126]{Ber2}, as well as in \cite[Ch.4, \S2, Prop.~9, p. 191]{MaG}.}

In general, using the subsheaf $\cI^1:=\cI$ we construct the following filtration of the structure sheaf:
\[
\widetilde\cO=\cI^0 \supset \cI^1 \supset \cI^2 \supset \ldots \supset \cI^m \supset \cI^{m+1}=0,
\]
and the associated  $\Zee$-graded sheaf of superalgebras
\[
\gr\widetilde\cO=\bigoplus\limits _{0\leq p\leq m}\gr\nolimits_p\widetilde\cO, \text{~~where $\gr_p\widetilde\cO=\cI^p/\cI^{p+1}$}.
\]
Here $\gr_0\widetilde\cO=\widetilde\cO_{\rd}$, and $\gr_1\widetilde\cO=\cF$ is a~locally free sheaf of $\widetilde\cO_{\rd}$-modules. Hence, it is clear that $\gr\widetilde\cO=\Lambda _{\widetilde\cO_{\rd}}(\cF)$. 

\hspace*{-1pt}Thus, to every supervariety $(M,\widetilde\cO)$ there is a~corresponding split supervariety $(M,\gr\widetilde\cO)$.

Denote by $\Der \widetilde\cO$ the sheaf of germs of derivations of the structure sheaf $\widetilde\cO$. Its stalk at point $w\in M$ is the Lie superalgebra $\Der _{\Cee}\widetilde\cO_w$. The sections of the sheaf $\Der \widetilde\cO$ are called \textit{vector fields} on the supermanifold $(M,\widetilde\cO)$. The space $\fd:=\Gamma(M,\Der \widetilde\cO)$ is naturally endowed with a~Lie superalgebra structure over $\Cee$.

If $(M,\widetilde\cO)$ is split, then $\Der \widetilde\cO$ is a~graded sheaf of Lie superalgebras whose homogeneous components are 
\[
\Der _p\widetilde\cO=\{\delta \in \Der \widetilde\cO \mid \delta \widetilde\cO_q \subset \widetilde\cO_{q+p}\quad \text{for all}\quad q\in \Zee\}.
\]
Therefore, $\fd$ is a~graded Lie superalgebra. In the general case we endow $\Der \widetilde\cO$ with the following filtration
\be\label{1}
\begin{array}{l}
\Der _{(-1)}\widetilde\cO=\Der \widetilde\cO,\\
\Der _{(p)}\widetilde\cO=\{\delta \in \Der \widetilde\cO \mid \delta \widetilde\cO \subset \cI^p,\ \ \delta \cI \subset \cI^{p+1}\}\quad \text{for all~~}p>0.
\end{array}
\ee

\ssbegin{Lemma}\label{l1} \textup{1)} The formulas ~\eqref{1} show that $\Der \widetilde\cO$ can be naturally considered as a~sheaf of filtered Lie superalgebras.

\textup{2)} The sheaves $\gr \Der \widetilde\cO$  and $\Der \gr\widetilde\cO$  are naturally isomorphic. \end{Lemma} 

\begin{proof} 1) It is easy to verify that 
\begin{equation}
\label{2}
\delta \cI^r \subset \cI^{r+p}\quad \text{for any $\delta \in \Der _{(p)}\widetilde\cO$ and $r\in \Zee$}.
\end{equation}
This implies that $[\Der _{(p)}\widetilde\cO,\Der _{(q)}\widetilde\cO]\subset \Der _{(p+q)}\widetilde\cO$.

2) To prove this, observe that every element $\delta \in \Der _{(p)}\widetilde\cO$ determines, thanks to eq.~\eqref{2}, linear mappings
\[
\tilde \delta _r:\cI^r/\cI^{r+1}\tto\cI^{r+p}/\cI^{r+p+1} \text{~~for all $r$}.
\]
It is subject to a~direct verification that the mappings $\tilde \delta _r$ form a~derivation $\tilde \delta \in \Der _p\gr\widetilde\cO$. The sheaf $\Der _{(p+1)}\widetilde\cO$ is the kernel of the mapping $\delta \mapsto \tilde \delta $. Thus we obtain an injective sheaf homomorphism 
\[
\gr_p \Der \widetilde\cO\tto   \Der_p \gr\widetilde\cO. 
\]
Using local splitness of supermanifolds, it is not difficult to show that this homomorphism is, moreover, surjective. Finally, a direct verification shows that this is an isomorphism of sheaves of graded Lie superalgebras.
\end{proof}

\sssbegin{Corollary}\label{C1} \textup{1)} The Lie superalgebra $\fd:=\Gamma (M,\Der \widetilde\cO)$ is filtered:
\[
\fd=\fd_{(-1)}\supset \fd_{(0)}\supset \fd_{(1)}\supset \ldots ,\text{~~where~~}\fd_{(p)}=\Gamma (M,\Der _{(p)}\widetilde\cO).
\]

\textup{2)} The following injective homomorphism of graded Lie superalgebras  is well-defined:
\[
{\gr\fd \tto \widetilde \fd:=\Gamma (M,\Der \gr\widetilde\cO)}.
\] 

If $M$ is compact, then $\dim _{\Cee}\fd<\infty $.\end{Corollary} 

\begin{proof} 1) The exact sequence
\[
0\tto \Der _{(p+1)}\widetilde\cO\tto \Der _{(p)}\widetilde\cO\tto \Der _p\gr\widetilde\cO\tto0,
\]
 existing thanks to Lemma~\ref{l1}, implies existence of the exact sequence 
\be\label{bez}
0\tto\fd_{(p+1)}\tto\fd_{(p)}\tto \Gamma (M,\Der_p\gr\widetilde\cO).
\ee

2) The sequence~\eqref{bez} implies  the existence of an injection 
\[
\gr\fd\tto \Gamma (M,\Der \gr\widetilde\cO).\hfill\qed
\]
\noqed\end{proof}

\section{Split supermanifolds}\label{S2}

Let  $\widetilde\cO=\Lambda(\cF)$, where  $\cF$ is a~locally free analytic sheaf of rank $m$ on an $n$-dimensional complex manifold $(M, \cO)$. Let $\cT$ denote the tangent sheaf  $\Der \cO$ on this manifold.

\ssbegin{Theorem}\label{Th1} \textup{1)} The sheaves  $\Der _p\widetilde\cO$ are locally free analytic. We have
\begin{equation}
\label{3}
\Der _{-1}\widetilde\cO \cong \Hom _{\cO}(\cF,\widetilde\cO)=\cF^*.
\end{equation}

\textup{2)} For $p\geq 0$, there is an exact sequence of sheaves
\begin{equation}
\label{4}
0\tto\cF^*\otimes \Lambda ^{p+1}(\cF)\tto \Der _p\widetilde\cO \overset \alpha \tto\cT \otimes \Lambda ^p(\cF)\tto 0.
\end{equation}
\end{Theorem}

\begin{proof} 1) It suffices to consider the case of the model supermanifold of dimension $n|m$. Clearly, in this case $\Der \widetilde\cO$  is a~free sheaf of $\widetilde\cO$-modules, a~basis of its sections consisting of
\[
\partial_{x_i}\quad (i=1,\ldots ,n),\qquad \partial_{\xi _j}\quad (j=1,\ldots ,m),
\]
where $x_1,\dots, x_n$ are coordinates in $\Cee^n$, and $\{\xi _1,\ldots ,\xi _m\}$ is any local basis of sections of the sheaf $\cF$. Therefore, $\Der _p\widetilde\cO$  is a~free sheaf of $\cO$-modules with a~basis of its sections being formed by 
\begin{gather*}
\xi _{i_1}\ldots \xi _{i_p}\partial_{x_i}\quad (i_1<\ldots <i_p,\ \ i=1,\ldots ,n)\\
\xi _{j_1}\ldots \xi _{j_{p+1}}\partial_{\xi _j}\quad (j_1<\ldots <j_{p+1},\ \ j=1,\ldots ,m)
\end{gather*}

Let $\delta\in\Der _p\widetilde\cO$. Since the sheaf $\widetilde\cO$ is generated by its subsheaves $\cO$ and $\cF$, the derivation $\delta$ is completely determined by its restrictions 
\[
\delta _0:=\alpha(\delta)=\delta|_\cO\text{~~and~~}\delta _1:=\beta(\delta)=\delta|_\cF.
\]
Obviously, 
\be\label{bn}
\delta _0\in\Hom_\Cee(\cO, \Lambda^p_{\cO}(\cF))\text{~~and~~}\delta _1\in\Hom_\Cee(\cF, \Lambda^{p+1}_{\cO}(\cF))
\ee
and
\begin{equation}
\label{5}
\begin{aligned}
\delta _0(\varphi \psi )&=(\delta _0 \varphi )\psi +\varphi (\delta _0 \psi )\\
\delta _1(\varphi s)&=(\delta _0 \varphi )s+\varphi (\delta _1s)\qquad \text{for any $\varphi ,
\psi \in \cO,s\in \cF$}.
\end{aligned}
\end{equation}

The other way around, for any pair $(\delta _0, \delta _1)$, see eq.~\eqref{bn}, which satisfies conditions ~\eqref{5}, there exists a~${\delta\in\Der _p\widetilde\cO}$ such that
\[
\alpha(\delta)=\delta _0\text{~~and~~}\beta(\delta)=\delta _1.
\]

For $p=-1$, we have $\alpha(\delta)=0$ and conditions ~\eqref{5} yields 
\[
{\beta(\delta)\in\Hom_\cO(\cF, \cO)=\cF^*}. 
\]
This implies the existence of the isomorphism ~\eqref{3}.

2) For $p\geq 0$, consider the sheaf homomorphism 
\[
\alpha: \Der _p\widetilde\cO\tto \Hom _{\cO}(\cO,\Lambda ^p_{\cO}(\cF)). 
\]
Obviously, $\Ker \alpha =(\Der _p\widetilde\cO)\cap \Der _{\cO}\widetilde\cO$, and therefore $\beta$ defines a~sheaf isomorphism
\[
\Ker\alpha\tto \Hom_{\cO}(\cF,\Lambda ^{p+1}_{\cO}(\cF))\cong \cF^*\otimes \Lambda ^{p+1}_{\cO}(\cF).
\]
To compute $\im \alpha $ we use an isomorphism 
\begin{equation}
\label{6}
\Hom _\Cee(\cO,\cG)\cong (\End _\Cee \cO) \otimes \cG
\end{equation}
for any locally free analytic sheaf $\cG$. Let $\{g_1,\ldots ,g_q\}$ be a~basis of sections of $\cG$ in a~neighborhood of a~ point $w\in M$. If $h\in \Hom _\Cee(\cO,\cG)$, then $h(\varphi )=\sum\limits _{1\leq i\leq q} \varphi _ig_i$, where $\varphi\in\cO_w$ and $\varphi_i\in\cO_w$  for all $i$. Setting $h_i(\varphi )=\varphi _i$ for all $i$, we see that $h_i\in (\End _\Cee \cO)_w$. The correspondence $h \mapsto \sum\limits _{1\leq i\leq q}h_i \otimes g_i$ does not depend on the basis $(g_i)_{i=1}^q$ and is the desired isomorphism ~\eqref{6}.

If $\cG=\Lambda ^p_{\cO}(\cF)$, then, for any given basis $\{g_i\}_{i=1}^q$, we can take local sections $\xi _{i_1} \dotsc \xi _{i_p}$ for $i_1 < \dotsb < i_p$. To the element $h\in \Hom _\Cee(\cO,\cG)$ the isomorphism ~\eqref{6} assigns 
\[
\text{$\sum\limits _{i_1<\ldots <i_p}h_{i_1 \ldots i_p}\otimes \xi _{i_1}\ldots \xi _{i_p}$,}
\]
 where $h_{i_1 \ldots i_p} \in \End _\Cee \cO$ is defined by the formula
\begin{equation}
\label{7}
h(\varphi )=\sum\limits _{i_1 <\ldots <i_p}h_{i_1 \ldots i_p}(\varphi )\xi _{i_1}\ldots \xi _{i_p}\quad \text{for any $\varphi \in \cO$}.
\end{equation}
Conditions ~\eqref{5}
 imply that if $h=\alpha (\delta )$, where $\delta \in \Der _p\cO$, then $h_{i_1 \ldots i_p} \in \Der \widetilde\cO$. 
 
 Thus, ${\im \alpha \subset \cT \otimes \Lambda ^p_\cO(\cF)}$, if we identify 
\[
\Hom _\Cee(\cO,\Lambda ^p_\cO(\cF))\simeq\End _\Cee \cO \otimes \Lambda ^p_\cO(\cF)\]
 by means of isomorphism~\eqref{6}. 

Conversely,  in formula~\eqref{7}, let $h\in \Hom _\Cee(\cO,\Lambda ^p_\cO(\cF))$  and ${h_{i_1 \ldots i_p} \in \cT}$. Then, in local coordinates $x_1,\ldots ,x_n$ on $M$, we have
\[
h_{i_1 \ldots i_p}=\sum\limits _{1\leq i\leq n} u_{i_1 \ldots i_p}^i \partial_{x_i},\quad \text{where}\quad u_{i_1 \ldots i_p}^i \in \cO.
\]
Clearly, 
\[
\alpha (\delta )=h\text{~~for~~} \delta =\sum\limits _{1\leq i\leq n} \ \ \sum\limits _{i_1<\ldots <i_p}\xi _{i_1}\ldots \xi _{i_p}\partial_{x_i} \in \Der _p\cO.
\]

Observe that for $p=0$ the sequence ~\eqref{4} is of the form
\begin{equation}
\label{8}
0\tto \End \cF\tto \Der _0\cO \overset \alpha \tto\cT\tto0.
\end{equation}
In particular, $\End \cF$ is a~sheaf of ideals in $\Der _0\cO$. The sequence of sheaves leads to the exact sequence of Lie algebras
\begin{equation}
\label{9}
0\tto \END F\tto\fd_0\tto \Gamma (M,\cT),
\end{equation}
where $\END F=\Gamma (M,\End \cF)$ is the Lie algebra of endomorphisms of the vector bundle $F$. 

Let $\varepsilon \in \fd_0$ be the element corresponding to the identity automorphism of $F$. Clearly, $\varepsilon $ is a~grading element of $\fd$, i.e., $[\varepsilon ,\delta ]=p \delta $ for any $\delta \in \fd_p$.
\end{proof}

\sssec{Example}\label{Ex2} Consider the split supermanifold $(M,\Omega )$, where 
\[
\Omega =\bigoplus\limits _{0\leq p\leq n}  \Omega ^p=\Lambda_{\cO}(\cT^*)
\]
is the sheaf of holomorphic forms on $(M,\cO)$. If $\cF=\cT^*$, then the sheaves which appear in Theorem~\ref{Th1} coincide with ${\Omega _\cT^p=\cT \otimes_{\cO} \Omega ^p} $ --- sheaves of vector-valued holomorphic $p$-forms, and we obtain the exact sequence which introduces maps $i_p$ and $\alpha _p$
\begin{equation}
\label{10}
0\tto \Omega _\cT^{p+1}\overset {i_p}\tto \Der _p \Omega \overset {\alpha _p}\tto \Omega _\cT^p\tto0.
\end{equation}

As was proved in \cite{2} (in the $C^\infty$ case, but this is inessential), this sequence splits. Therefore, we obtain the following theorem which I prove for completeness.

\ssbegin{Theorem}\label{Th2} For any $p \geqslant -1$, we have
\[
\Der _p \Omega \cong \Omega _\cT^p \oplus \Omega _\cT^{p+1}.
\]
\end{Theorem}

\begin{proof} Let $\omega \in \Omega _\cT^p$. Consider $\delta _\omega :=i_{p-1}(\omega ) \in \Ker \alpha _{p-1}$, and set $\delta :=[\delta _\omega ,d]$, where $d:\Omega \tto \Omega $ is the exterior differential, a~section of the sheaf $\Der _1 \Omega $. Then, $\alpha _p(\delta )=\omega $.

Indeed, let $\omega =\sum\limits _i \delta _i \omega _i$, where $\delta _i\in \cT$ and $\omega _i\in \Omega ^p$. Then, for any $\varphi \in \cO$, we have
\[
\delta \varphi =\delta _\omega (d \varphi )=\sum\limits _i(d \varphi )(\delta _i)\otimes \omega _i=\sum\limits _i \delta _i(\varphi )\otimes \omega _i=\omega (\varphi ).
\]
Therefore, the mapping $\omega \mapsto \delta =[\delta _\omega ,d]$ splits the sequence ~\eqref{10}. It is easy to see that the subsheaf of $\Der _p \Omega $ defined by the splitting, which is complementary to $\Omega _\cT^{p+1}$, coincides with the centralizer of~$d$. (The first description of  the centralizer of~$d$ in terms of Lie superalgebras and invariant differential operators is due to Grozman; for details, see the arXiv version of \cite{G} reproduced in this Special Volume \textit{D.L.}).
\end{proof}

\ssbegin{Corollary}\label{Cor} There exist isomorphisms 
\[
\fd=\Gamma (M,\Der _p \Omega )\cong \Gamma (M,\Omega _\cT^p)\oplus \Gamma (M,\Omega _\cT^{p+1}).
\]
In particular, $\fd_{-1}\cong \Gamma (M,\cT)$. 
\end{Corollary}

The latter isomorphism is described as follows: every vector field $v\in \Gamma (M,\cT)$ determines a~derivation  $i_v \in \fd_{-1}$, called \textit{inner derivation by} (or a~\textit{convolution with}) $v$. Further, we have a~decomposition into a~semidirect sum of Lie algebras
\[
\fd_0 \cong \Gamma (M,\cT)\oplus \END T,
\]
where $\END T$ is an ideal, and the Lie subalgebra  $\Gamma (M,\cT)$  is embedded into $\fd_0$ by means of injective homomorphism $v \mapsto \theta_v$, where $\theta_v$ is the \textit{Lie derivative along the field} $v$. Finally, 
\[
\fd_1 \cong \END T \oplus \Gamma (M,\Omega _\cT^2).
\]
Under this isomorphism the identity automorphism $\varepsilon \in \END T$ corresponds to the exterior differential $d\in \fd_{1}$.

Consider the graded subspace $\hat \fd=\hat \fd_{-1}\oplus \hat \fd_0 \oplus \hat \fd_1 \subset \fd$, where
\[
\begin{array}{l}
\hat \fd_{-1}=i_{\Gamma (M,\cT)}=\fd_{-1},\ \
\hat \fd_0=\theta_{\Gamma (M,\cT)}\oplus \Cee \varepsilon \subset \fd_0,\ \
\hat \fd_1=\Cee d\subset \fd_1.
\end{array}
\]
The classical relations (here $p(i_v)=\od$, i.e., $i_v$ is odd for all vector fields $v$)
\begin{gather*}
[i_v,i_w]=0,\quad [\theta_v,\theta_w]=\theta_{[v,w]},\quad [\theta_v,i_w]=i_{[v,w]},\\
[d,i_v]=\theta_v,\quad [d,\theta_v]=0 \quad \text{for any $v,w\in \Gamma (M,\cT)$},\\
[d,d]=0,\quad [\varepsilon ,\delta ]=p \delta \quad \text{for any $\delta \in \fd_p$}
\end{gather*}
immediately imply that $\hat \fd$ is a~Lie subsuperalgebra in $\fd$.

\section{Transitive Lie superalgebras}\label{S3}

Let us deduce a sufficient condition for the Lie superalgebra of vector fields $\Gamma (M,\Der \widetilde \cO)$ on the supermanifold $(M,\widetilde \cO)$ to be transitive. As before let $\cJ$ be the subsheaf of ideals generated by odd elements. The locally free sheaf $\cF^*=(\cJ/\cJ^2)^*$ on the complex manifold $(M,\widetilde \cO_{\rd})$ will be called \textit{odd tangent sheaf} and denoted $\cT_{\overline 1}$. Let $T_{\overline 1}$ be the corresponding holomorphic vector bundle over $(M,\widetilde \cO_{\rd})$ called the \textit{odd tangent bundle}.

\ssec{Split case, i.e., $\widetilde \cO=\Lambda(\cF)$} By Theorem \ref{Th1}, $\fd_{-1}\cong \Gamma (M,\cT_{\overline 1})$.

\sssbegin{Lemma}\label{L2} Let $(M,\widetilde \cO)$ be a~split supermanifold. Suppose  a~global holomorphic section passes through every point of any fiber of the bundle $T_{\overline 1}$. 

Then, the graded Lie superalgebra $\fd=\Gamma (M,\Der \widetilde \cO)$ satisfies  condition ~\eqref{12} for all $p>0$. 

If this condition holds for $p=0$ as well, i.e., if the adjoint representation of $\fd_0$ in $\fd_{-1}$ is exact, then $\fd$ is transitive. \end{Lemma}

\begin{proof} Let $\delta \in \fd_p$ for $p>0$ and $\gamma \in \fd_{-1}$. Then, 
\begin{gather}
\label{13}
[\gamma ,\delta ](\varphi) =\gamma (\delta \varphi )\text{~~for any $\varphi \in \widetilde \cO$},\\
[\gamma ,\delta ](s)=\gamma (\delta s)+(-1)^{p+1}\delta (\gamma s)\text{~~for any $s\in \cF$}.\hfill \qed
\label{14}
\end{gather}
\noqed\end{proof}

If $[\gamma ,\delta ]=0$ for all $\gamma \in \fd_{-1}$, then $\gamma (\delta \varphi )=0$ for all $\varphi \in \cO$. This follows from ~\eqref{13}. Since $\delta \varphi \in \Lambda ^p(\cF)$, then under Lemma's hypothesis for $p>0$, we have $\delta \varphi =0$  for all $\varphi \in \cO$. But then Eq.~\eqref{14} implies that $\gamma (\delta s)=0$ for all $\gamma \in \fd_{-1}$  and $s\in \cF$, so $\delta s=0$ for all $s\in \cF$; i.e., $\delta =0$.

Observe that, in general, the Lie superalgebra $\fd$ does not satisfy condition ~\eqref{12} for $p=0$. Under assumptions of the Lemma we can only claim that the action of the ideal $\END F=\END T_{\overline 1}$, see ~\eqref{9}, on $\fd_{-1}=\Gamma (M,\cT_{\overline 1})$ is exact.

\ssec{General case} Let $\fd=\Gamma (M,\widetilde \cO)$. Consider the graded Lie superalgebra 
\[
\tilde \fd:=\Gamma (M,\Der \gr\widetilde \cO)\simeq\Gamma (M,\gr \Der \widetilde \cO).
\]
By item 2) of Corollary \ref{C1},  there is an injective homomorpshism $\gr\fd\tto \tilde \fd$. 

In particular, ${\gr_{-1}\fd=\fd/\fd_{(0)}}$ is identified with %the subspace 
$\tilde \fd_{-1}=\Gamma (M,\cT_{\overline 1})$.

\sssbegin{Lemma}\label{L3}  For the Lie superalgebra $\fd:=\Gamma (M,\Der\widetilde \cO)$ to be transitive the following two conditions are sufficient:

\textup{a)} Given any point of any fiber of the bundle $T_{\overline 1}$, there is a~section of the subspace 
\[
{\fd_{-1}=\fd/\fd_{(0)}\subset \Gamma (M,\cT_{\overline 1})}
\]
whose image contains it. 

\textup{b)} The adjoint action of the Lie algebra $\fd_0:=\fd_{(0)}/\fd_{(1)}$  in $\fd_{-1}$  is exact. \end{Lemma}

\begin{proof} By repeating almost verbatim the proof of Lemma \ref{L2} we see that the subalgebra $\gr\fd$ of the Lie superalgebra $\tilde \fd$ satisfies conditions ~\eqref{12} for any $p>0$. The case $p=0$ is handled by condition of item b). \end{proof}

\section{Homogeneous vector bundles}

Let $p:E\tto M$  be a~holomorphic vector bundle over a~complex manifold $(M,\cO)$. A fiber-wise linear biholomorphic mapping $E\tto E$ is said to be an \textit{automorphism} of the bundle $E$. Let $A(E)$  be the group of all automorphisms of $E$. Obviously, every automorphism $a$ of the bundle $E$ determines an automorphism $p(a)$ of the base $M$ of $E$. We obtain an exact sequence of groups
\begin{equation}
\label{15}
e\tto \Aut E\tto A(E)\overset p\tto \Aut M,
\end{equation}
where $\Aut E \subset A(E)$ is the normal subgroup consisting of the automorphisms sending every fiber into itself. If $M$ is compact, then the sequence ~\eqref{15} consists of complex Lie groups and their homomorphisms, see \cite{Mo}.

Observe that it is possible to describe the automorphisms of the bundle $E$ in terms of the corresponding sheaf $\cE$. Namely, every $a\in A(E)$ determines an automorphism $\tilde a$ of the sheaf $\cE$ over the automorphism $p(a)$ of the base $M$, i.e., determines an isomorphism of sheaves $\tilde a:\cE\tto p(a)^*\cE$. The latter isomorphism is given on local sections $s$ by the formula
\begin{equation}
\label{16}
\tilde a(s)(w)=a^{-1}(s(p(a)(w)))\text{~~for any $w\in M$}.
\end{equation}

Consider the sheaf $\cE^*=\Hom _{\cO}(\cE, \cO)$ corresponding to the dual bundle $E^*$. The space of sections $\Gamma (U,\cE^*)$ over any open set $U \subset M$ can be identified with the subspace in $\Gamma (p^{-1}(U),\cO_E)$ consisting of the functions linear on fibers. This gives an embedding $\cE^* \subset p_*\cO_E$.

A vector field on $E$ with a~projection to $M$ and sending $\cE^*$ into itself will be called an \textit{infinitesimal automorphism of the bundle} $E$. The infinitesimal automorphisms determine a~sheaf of complex Lie algebras $\cA(E)$ on $M$. Projection to the base yields a~sheaf homomorphism $\pi :\cA(E)\tto\cT$, where $\cT$ is the tangent sheaf on $M$. If $M$ is compact, then the Lie algebra $\mkA (E)=\Gamma (M,\cA(E))$ is tangent to the Lie group $A(E)$, and the homomorphism ${\pi :\mkA (E)\tto \Gamma (M,\cT)}$ coincides with $dp$.

\ssbegin{Lemma}\label{L4}  The sheaf $\cA(E)$ is naturally isomorphic to the sheaf $\Der _0\widetilde \cO$ of degree $0$ derivations of ${\widetilde \cO:=\Lambda(\cE^*)}$.   Under this isomorphism $\pi $ turns into the homomorphism $\alpha $ from sequence~\eqref{8}. \end{Lemma}

\begin{proof}Every $\delta \in \cA(E)$ determines a~pair
\[
\text{$\delta _0=\pi (\delta ) \in \cT$ and $\delta _1=\delta |_{\cE^*} \in \End _\Cee\cE^*$,}
\]
 so that conditions~\eqref{5} are satisfied. Let $\hat \delta \in \Der _0\widetilde \cO$ be the derivation corresponding to the pair $(\delta _0,\delta _1)$, see proof of Theorem \ref{Th1}. It is easy to see that the mapping $\delta \mapsto \hat \delta $ is an injective homomorphism of sheaves of Lie algebras. 

To prove surjectivity of the mapping $\delta \mapsto \hat \delta $, consider $\hat \delta \in (\Der _0\widetilde \cO)_w$ determined by a~pair $(\delta _0,\delta _1)$ satisfying conditions~\eqref{5}. Let $U \subset M$ be an open neighborhood of any point $w\in M$ over which $E$ admits trivialization $p^{-1}(M)=U \times \Cee^m$. Set $\delta \varphi =\delta _0 \varphi $, where $\varphi \in \cO |_U$, and $\delta  \ell =\delta _1 \ell $, where $\ell \in (\Cee^m)^*$.These conditions completely determine an element $\delta \in \cA(E)_w$ which turns into $\hat \delta $ under this mapping~$\widehat{\ }$.
\end{proof}

Now let $(M,\cO)$ be a~homogeneous space of a~connected complex Lie group $G$, i.e., let there be given a~transitive analytic $G$-action on $M$. Fix a~point  $w_0\in M$, and let $P$ be the stationary subgroup of $w_0$ in $G$. 

As is well-known, $M$ can be naturally identified with the manifold of the left cosets $G/P$. Let $M$ be compact.

Recall that a~holomorphic vector bundle $p:E\to M$ is \textit{homogeneous} (under a~$G$-action) if there is an analytic $G$-action  by automorphisms of the bundle $E$ whose projection is a~given $G$-action on $M$. Equivalently, there exists an analytic homomorphism $\hat t:G\to A(E)$ such that $p \circ \hat t=t $ is a~homomorphism $G\to \Aut\ M$ which determines the given $G$-action on~$M$.

Associated to a~homogeneous bundle $E$ there is a~linear representation $\varphi :g \mapsto \hat t(g) |_{E_{w_0}}$ of $P$ on the space $E_{w_0}$. This representation completely determines the bundle $E$ together with the {$G$-action} on~$E$. For the role of $\varphi $ one can take any finite-dimensional linear analytic representation of $P$. We write $E_\varphi $ to denote the homogeneous vector bundle over $M$ corresponding to the representation $\varphi $. The corresponding locally free sheaves $\cE_\varphi $ are also called \textit{homogeneous}.

The $G$-action on $E_\varphi $ determines an analytic linear representation $\Phi \!:\! G\to GL(\Gamma (M,\cE_\varphi ))$ given by the formula
\[
(\Phi (g)s)(w)=gs(g^{-1}w)\text{~~ for any $g\in G,\quad w\in M$}.
\]
The representation $\Phi $ of $G$ is called \textit{induced} by the representation $\varphi $ of $P$.

\sssbegin{Example}\label{ex} The tangent bundle $T$ over $M$ is endowed with a~natural $G$-action, and therefore is homogeneous. It corresponds to the linear representation $\tau :P\to GL(T_{w_0}(M))$ given by the formula 
\[
\tau (g)=dt(g)_{w_0} \text{~~ for any $g\in P$}.
\]
\end{Example}

In what follows we consider the case where $M$ is the homogeneous flag variety, i.e., where $G$ is a~connected semisimple (or reductive) complex Lie group and $P$ is its parabolic subgroup. In this case, the space of sections of the homogeneous vector bundle and the induced representation are described by the famous Bott--Borel--Weil theorem, see \cite{B}. Let us recall it.

Let $T$ be a~\textit{maximal torus} in a~\textit{Borel subgroup} $B$ of $G$ and  let $\ft$ and $\fb$ be the respective Lie algebras. Let $R$ be the \textit{root system} with respect to $T$ (or $\ft$); let $R^+$ be the \textit{set of positive roots} corresponding to $\fb$ and $\Pi$ the \textit{set of simple roots} and let 
\[
\left\{h_\alpha:=\frac{2 \alpha}{( \alpha, \alpha)} \mid \alpha \in R\right\}\subset \ft 
\]
be the set of \textit{coroots}, i.e., the system of roots dual to $R$. The weight $\lambda \in \cI^*$ is called \textit{dominant} if $\lambda (h_\alpha )\geqslant 0$ for all $\alpha \in \Pi$.

\ssbegin{Theorem}[Bott--Borel--Weil]\label{BWB}\label{Th3} Let $P$ be a~parabolic subgroup of $G$ containing the opposite to $B$ Borel group $B^-$, and let $M=G/P$. Let $\varphi $ be an irreducible analytic representation of $P$ with highest weight~$\lambda $. 

Then, $\Gamma (M,\cE_\varphi )\ne 0$ if and only if $\lambda $ is dominant. In this case, the induced representation $\Phi $ in $\Gamma (M,\cE_\varphi )$ is irreducible with highest weight~$\lambda$. \end{Theorem}

This theorem is applicable only if $\varphi $ is irreducible or completely reducible. However, one often encounters homogeneous bundles arising from representations which are not completely reducible. The following Lemma is useful for studying them.

\ssbegin{Lemma}\label{L5} Let $M=G/P$ be a~homogeneous flag variety, $\varphi $ a~holomorphic linear representation of $P$. If the induced representation $\Phi $ of $G$ contains an irreducible representation with highest weight $\lambda $ of multiplicity $t$, then $\lambda $ is contained in the set of weights of $\varphi$ of multiplicity $\geqslant t$. \end{Lemma}

\begin{proof} Let us perform induction by the length $\ell (\varphi )$ of the Jordan-Hoelder series of  $\varphi $. If $\ell (\varphi )=1$, then $\lambda $ is the highest weight of $\varphi $ by Theorem \ref{Th3}. Let the Lemma be proved for representations $\varphi '$ such that $\ell (\varphi ')<\ell (\varphi )$. Assuming that $\varphi $ is reducible, consider a~proper subrepresentation $\varphi _1$ of $\varphi $ and the quotient representation $\varphi _2$. We obtain an exact sequence of sheaves
\[
0 \tto \cE_{\varphi _1} \tto \cE_\varphi \tto \cE_{\varphi _2} \tto 0,
\]
which defines the exact sequence of $G$-modules 
\[
0\tto \Gamma (M,\cE_{\varphi _1}) \tto \Gamma (M,\cE_\varphi ) \tto \Gamma (M,\cE_{\varphi _2}).
\]
Since every finite-dimensional analytic representation of $G$ is completely reducible, the irreducible representation with highest weight $\lambda $ is contained in $\Gamma (M,\cE_{\varphi _1})$ and $\Gamma (M,\cE_{\varphi _2})$ with total multiplicity $\geqslant t$. By induction hypothesis the statement of the Lemma follows because the set of weights of $\varphi $ is the union of the sets of weights of $\varphi _1$ and $\varphi _2$.
\end{proof}

\sssbegin{Corollary}\label{col}  If $\varphi $ has no dominant weights, then the induced representation $\Phi $ in $\Gamma (M,\cE_\varphi )$ is trivial. Moreover, ${\Gamma (M,\cE_\varphi )=0}$. \end{Corollary}

\section{Split homogeneous supervarieties}\label{Ssph}

In this section we consider a~particular class of split supervarieties associated with homogeneous flag varieties. Let $G$ be a~connected semisimple (or reductive) complex Lie group, $P$ a~ complex Lie subgroup, and $E_\varphi $ a~ homogeneous vector bundle over $M=G/P$ determined by a~representation $\varphi $ of $P$. Consider a~split supermanifold $(M,\widetilde \cO)$, where $\widetilde \cO=\Lambda(\cE_\varphi)$. The sheaf $\widetilde \cO$ corresponds to the homogeneous vector bundle $\Lambda (E_\varphi) =E_{\Lambda(\varphi)}$, i.e., $\widetilde \cO=\cE_{\Lambda (\varphi)}$. 

The sheaf $\Der \widetilde \cO$ is also homogeneous. Indeed, every element $g\in G$ corresponds to an automorphism $\widetilde {\hat t(g)}$ of $\widetilde \cO$ over an automorphism $t(g) \in \Aut M$, see eq.~\eqref{16}. By setting
\begin{equation}
\label{17}
g \delta =\widetilde {\hat t(g)} \circ \delta \circ \widetilde {\hat t(g)}{}^{-1}\text{~~ for any $\delta \in \Der \widetilde \cO$}
\end{equation}
we obtain the desired $G$-action on the locally free sheaf $\Der \widetilde \cO$.  Also notice that the sheaves in the exact sequence~\eqref{4} correspond to the homogeneous vector bundles $E_\varphi^* \otimes \Lambda ^{p+1}(E_\varphi)$ and $T \otimes \Lambda ^p(E_\varphi)$, while the homomorphisms in this sequence are $G$-equivariant.

By Lemma \ref{L4} we can identify the Lie algebra $\fd_0=\Gamma (M,\Der _0\widetilde \cO)$ with $\mkA (E_\varphi ^*)$. Since the bundle $E_\varphi ^*$ is also homogeneous, we have a~Lie algebra homomorphism $d \hat t:\fg \tto \fd_0$. On the other hand, the $G$-action ~\eqref{17} on $\widetilde \cO$ determines a~representation $\Psi :G\tto GL(\fd)$.

\ssbegin{Lemma}\label{usfL}\label{L6} Let $u\in \fg$ and $\psi =d \Psi :\fg\tto \fgl (\fd)$. 
Then,  $[d \hat t(u),\delta ]=\psi (u)\delta $. \end{Lemma}

\begin{proof} This follows directly from ~\eqref{17}: instead of $g$,  substitute in eq.~\eqref{17}  the curve $g(t)\in G$ with tangent vector $u$ at $t=0$ passing through $e:=g(0)$, and differentiate both parts with respect to $t$ at $t=0$. \end{proof} 

A split supermanifold $(M,\Lambda(\cE_\varphi ))$ will be called \textit{homogeneous} if ${T_{\overline 1}=E_\varphi ^*=E_{\varphi ^*}}$ has ``sufficiently many'' holomorphic sections. This means that 
every point of every fiber of this bundle is contained in the image of some section. 

In what follows $G$ is semisimple and $P$ is its parabolic subgroup.

\ssbegin{Lemma}\label{L7} Let  $\varphi $ be irreducible and  $\lambda $ the highest weight of representation $\varphi ^*$. The following properties are equivalent:
\begin{equation*}
\begin{array}{ll}
i)&(M,\Lambda (\cE_\varphi ))\text{~~is homogeneous};\\
ii)&\lambda\text{~~is dominant};\\
iii)&\fd_{-1}=\Gamma (M,\cE_{\varphi ^*})\ne 0. \\
\end{array}
\end{equation*}
Under these conditions the induced representation $\Psi _{-1}$ of $G$ in $\fd_{-1}$ is irreducible. \end{Lemma}

\begin{proof} By Theorem~\ref{Th3}, it remains to prove that if $\lambda$ is dominant, then $T_{\overline 1}$ has sufficiently many holomorphic sections. Obviously, the restriction map $r_{w_0}:\Gamma (M,\cE_{\varphi ^*}) \to (E_{\varphi ^*})_{w_0}$ is a~$P$-module homomorphism. Therefore, either $r_{w_0}$ is surjective or $r_{w_0}=0$. But for any $g\in G$ and $s\in \Gamma (M,\cE_{\varphi ^*})$, we have
\begin{equation}
\label{18}
r_{gw_0}(s)=gr_{w_0}(\Phi (g^{-1})s),
\end{equation}
where   $\Phi$ is the representation induced by $\varphi ^*$. Therefore, $r{w_0}=0$ implies $r_w=0$ for all $w\in M$, which is impossible. Thus $r_w$ is surjective for all $w\in M$.
\end{proof}
 
Consider the following commutative diagram whose upper line is the exact sequence ~\eqref{9} for the bundle $F=E_\varphi ^*$:
\begin{equation}
0\tto \END E_\varphi \tto \begin{array}[t]{@{}rcl@{}}
\fd&\overset \alpha \tto &\Gamma (M,\cT)\\[-2pt]
\raisebox{-5pt}{$\scriptstyle d \hat t$}\!\!\nwarrow \hspace*{-8pt}&&\hspace*{-8pt}\nearrow \!\!\raisebox{-5pt}{$\scriptstyle dt$}\\[-7pt]
&\fg
\end{array}
\label{19}
\end{equation}

\ssbegin{Theorem}\label{Th4} Let $M=G/P$ be a~flag variety of a~ semisimple Lie group $G$. Let $\varphi $ be an irreducible finite-dimensional analytic linear representation of $P$; let $\widetilde \cO=\Lambda (\cE_\varphi)$ and ${\alpha (\fd_0)=dt(\fg)}$, where $\alpha$ is defined in eq.~\eqref{19}. 

Let $G=G_1 \ldots G_r$ be a~decomposition of $G$ into simple factors, and let $\Pi =\Pi_1 \sqcup \dotsm \sqcup \Pi_r$ be the corresponding decomposition of the system of simple roots. Let the highest weight $\lambda $ of the representation $\varphi ^*$ be dominant, and for every $i$ assume that there exists a~$\beta \in \Pi _i$ such that $\lambda (h_\beta )>0$. 

Then, the graded Lie superalgebra $\fd=\Gamma (M,\Der \cO)$ is transitive and irreducible. \end{Theorem}

\begin{proof} By Lemma \ref{L7} and  since $\lambda $ is dominant,  the condition of Lemma \ref{L2} is satisfied. Therefore,  it remains to show that the adjoint representation of $\fd_0$ in $\fd_{-1}=\Gamma (M,\cE_{\varphi ^*})$ is irreducible and exact. 

The irreducibility easily follows from Lemmas \ref{L6} and \ref{L7}.

Observe that the homomorphism $d \hat t$ is injective. Indeed, otherwise $\fg_i \subset \Ker d \hat t$ for some~$i$. Hence, $dt(\fg_i)=0$ which implies that $G_i \subset P$. Furthermore, $\fg_i \subset \Ker \psi $ by Lemma \ref{L6}. 

Since $r_{w_0}:\Gamma (M,\cE_{\varphi ^*}) \to (E_{\varphi ^*})_{w_0}$ is surjective (Lemma \ref{L7}), it follows that $G_i \subset \Ker \varphi^*$, which contradicts the hypothesis.
\end{proof}

Let us identify $\fg$ with $d \hat t(\fg)$ by means of $d \hat t$. Form eq. ~\eqref{19} we see that $dt=\alpha |_\fg$ and ${\alpha (\fd_0)=\alpha (\fg)}$.

Consider the ideal $\fg_0=\fg \cap \END E_\varphi $ of $\fg$. We see that $\fg=\fg_0 \oplus \tilde {\fg}$, where $\tilde {\fg}$ is also an ideal. Obviously, our assumptions imply that there is a~decomposition into a~semi-direct sum
\[
\fd_0=\END E_\varphi \ltimes \tilde \fg,
\]
where $\END E_\varphi  $ is an ideal and $\tilde \fg$ a~subalgebra. Let us show that this sum is actually direct. 

Consider the induced representation $\psi $ of the Lie algebra $\fg$ in $\END E_\varphi $. The weights of the induced representation $\varphi \otimes \varphi ^*$ are of the form $\mu =\nu _1-\nu _2$, where $\nu _1$, $\nu _2$ are weights of $\varphi $. Since $\varphi $ is irreducible, $\mu $ can be expressed in terms of the set $\Pi _P \subset \Pi $ corresponding to the semisimple part of $P$. But the roots of any subset $\Pi _i$ should enter the expression of any dominant weight in terms of $\Pi $ with either positive or zero coefficients, see \cite[\S13, Exc. 8]{3} and \cite{9}. 

Since $(dt)|_{\tilde \fg}$ is injective, $G_i\not \subset P$ for all $i$ such that $\fg_i \subset \tilde \fg$. Therefore, if $\mu $ is dominant, then  $\mu (h_\beta )=0$ for all $\beta \in \Pi _i$ and any such $i$.

Since $\varphi \otimes \varphi ^*$ is completely reducible, we can apply Theorem \ref{Th3} which implies that $\psi (\tilde \fg)=0$, and $[\tilde \fg,\END E_\varphi ]=0$ thanks to Lemma \ref{L6}.

The homogeneity of the supermanifold $(M,\widetilde\cO)$ proved above implies that the action of $\END E_\varphi =\END E_{\varphi ^*}$ in $\fd_{-1}$ is exact. The radical $\fr$ of the Lie algebra $\fd_0$ is contained in $\END E_\varphi $ and is non-zero since $\Cee \varepsilon \subset \fr$, see \S~2. Further, the Lie subalgebra $\ad_{\fd_0} \subset \fgl(\fd_{-1})$ is irreducible, and the dimension of its radical should be at most $1$. Hence, $\fr=\Cee \varepsilon $.

Therefore, $\END E_\varphi =\Cee \varepsilon \oplus \fh$  and $\fd_0=\END E_\varphi \oplus \tilde \fg$ are reductive Lie algebras with 1-dimensional centers. Any ideal of $\fd_0$ is contained either in $\END E_\varphi $ or in $\tilde \fg$. Thanks to our assumption on $\lambda $ and Theorem \ref{Th3}, the action of the Lie algebra $\tilde \fg$ in $\fd_{-1}$ is exact, and so is the $\fd_0$-action.

Observe that the assumptions of Theorem \ref{Th4} are satisfied, e.g., if $G$ is simple and coincides with $\Aut ^\circ M$, where ${}^\circ$ singles out the connected component of the unit, and $\varphi $ is a~non-trivial irreducible representation such that the highest weight of $\varphi ^*$ is dominant.

On the other hand, the version of the theorem proved above, which does not presuppose that $G$-action on $M$ is faithful, is useful in applications. As is well-known, the cases where the inclusions between Lie algebras $\alpha (\fd_0)$, $(dt)(\fg)$ and $\Gamma (M,\cT)$ are strict do occur  very seldom.

Observe that, if $\varphi $ is irreducible and the $G$-action on $M$ is locally faithful, then we have ${\END E_\varphi =\Cee \varepsilon}$, see \cite{9}.

In conclusion, let us compute the Lie superalgebras $\fd$ for the supermanifolds $(M,\Omega )$ of Examples~ \ref{Ex2} and \ref{ex}, where $M$ is a~flag variety. Obviously, such a~supermanifold is homogeneous for any compact complex homogeneous manifold $M$, although the representation $\varphi $ that determines $(M, \Omega)$ is seldom completely reducible.

\ssbegin{Theorem}\label{Th5} Let $M=G/P$ be a~flag variety of the simple Lie group $G$. 
Then, $\END T=\Cee \varepsilon $ and $\Gamma (M,\Omega _\cT^p)=0$ for $p>1$; the Lie superalgebra $\fd=\Gamma (M,\Der \Omega )$  coincides with its subsuperalgebra $\hat \fd$, see $\S2$. \end{Theorem}

\begin{proof} The equality $\END T=\Cee \varepsilon $ is proved in \cite{4}. It can also be deduced from Theorem \ref{Th3} by using certain facts of Lie algebra theory. 

Let $p>1$. Then, the weights of the representation $\tau \otimes \Lambda ^p (\tau ^*)$ which determines the sheaf $\Omega _\cT^p$ are of the form 
\[
\alpha -\beta _{i_1}-\ldots -\beta _{i_p},\text{~~where~~}\alpha, \beta _{i_1},\ldots ,\beta _{i_p} \in \Pi^+
\]
Obviously, such a~weight can not be dominant. Applying item 1) of Corollary \ref{C1}, we see that $\Gamma (M,\Omega _\cT^p)=0$.

The Lie superalgebra $\fd$ spoken about in Theorem \ref{Th5} is isomorphic to $\fvect(0|m)$.
\end{proof}

\end{paper}